# DOES THE ISHANGO BONE INDICATE KNOWLEDGE OF THE BASE 12? AN INTERPRETATION OF A PREHISTORIC DISCOVERY, THE FIRST MATHEMATICAL TOOL OF HUMANKIND

Vladimir PLETSER
European Space Research and Technology Centre,
European Space Agency, P. O. Box 2200 Noordwijk, The Netherlands

**Abstract**

In the early fifties, the Belgian Prof. J. de Heinzelin discovered a bone in the region of a fishermen village called Ishango, at one of the farthest sources of the Nile, on the border of Congo and Uganda. The Heinzelin's Ishango bone has notches that seem to form patterns, making it the first tool on which some logic reasoning seems to have been done. In this paper a new interpretation is proposed for these patterned notches, based on a detailed observation of their structure. It can be called the "slide rule"-reading, in contrast to former "arithmetic game" and "calendar" explanations. Additional circumstantial evidences are given to support the hypothesis that the Ishango bone is a primitive mathematical tool using the base 12 and sub-bases 3 and 4.

## 1. Introduction.

The Ishango bone is a 10-cm long curved bone, first described by its discoverer, Prof. J. de Heinzelin [deH1, deH2]. He found the tiny bone about fifty years ago, among harpoon heads at a certain depth in stratified sand in the area of the Ishango fishermen village on the shores of the Semliki river, not far from the present border between Congo and Uganda. The bone has encased in its narrower (top) extremity a fragment of quartz, protruding by 2 mm, most probably for tattooing or engraving purposes. Dating of the bone is somehow difficult and most likely indicate an age of 20 000 years (although other indications point toward an age of 90 000 years), which makes it the first mathematical tool of Mankind.

Most interestingly, the bone carries 167 or 168 notches distributed in three columns along the bone length. One of the columns is typically central along the most curved side of the bone. The columns are called M, G and D following the initial of the French words Middle (*Milieu*), Left (*Gauche*), and Right (*Droite*). Within each column, the notches are grouped like depicted in Figure 1, most likely not at random. The M column shows from top to bottom eight groups of respectively 3, 6, 4, 8, 9 or 10, 5, 5, and 7 notches. The G and D columns show each four groups respectively of 11, 13, 17, 19 and of 11, 21, 19, 9 notches. The notches are approximately parallel within each group, sometimes of different lengths and of different orientations. Some uncertainty remains on the number of notches for the Me(10) and Mf(5) groups of 10 and 5 notches, as part of the bone was damaged earlier or by infiltrating rainwater. The last bottom notch of Me(10) group is somehow separated from the rest of the group and is interrupted in its middle. The first top notch of the Mf(5) group is not clearly visible when compared to the other 4 notches in the group.

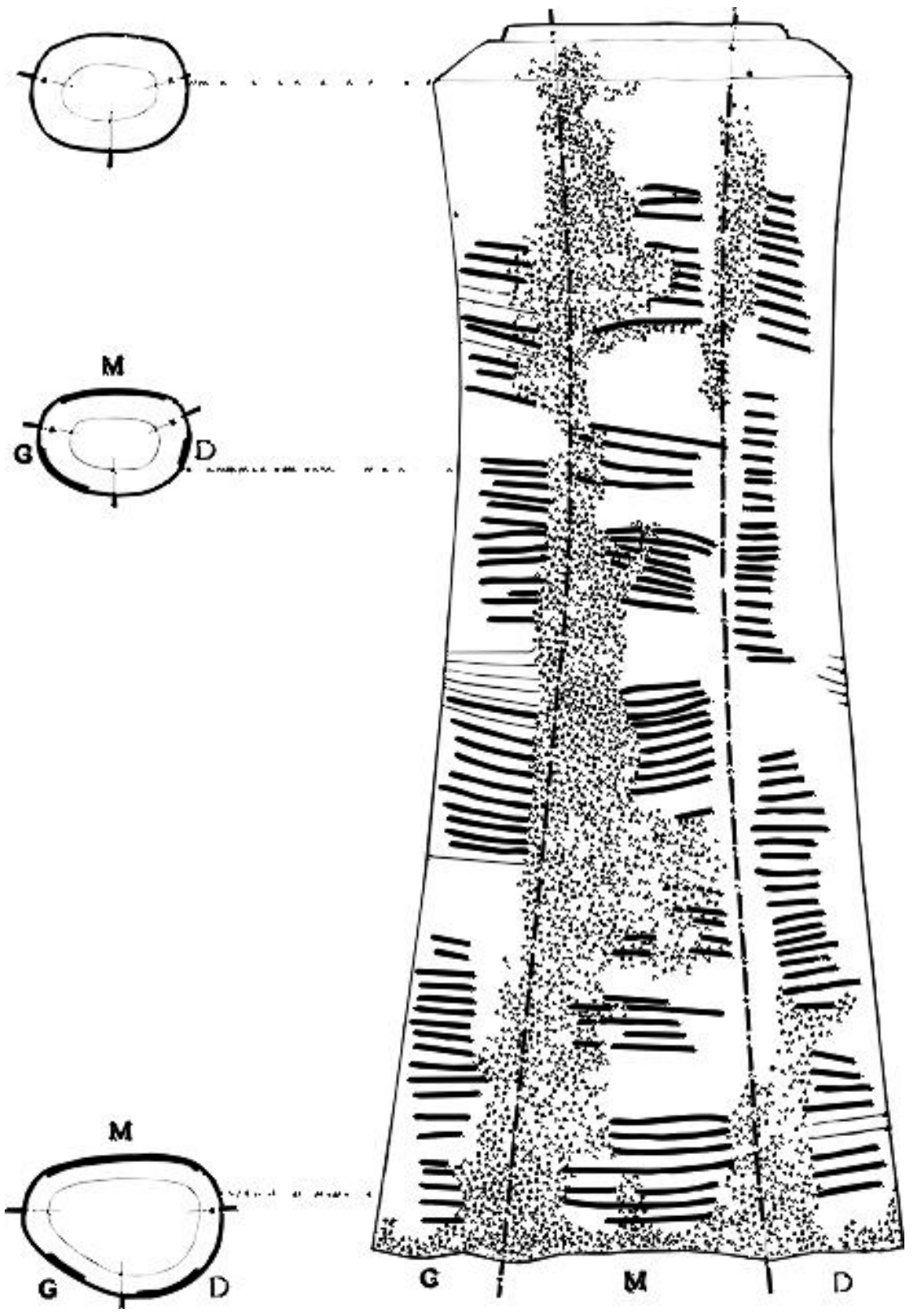

Figure 1: De Heinzelin's faithful detailed drawing of the Ishango bone.

Nevertheless, replacing the group of notches by the number of notches yield a table of simple numbers like shown in Figure 2. The M column includes numbers less than or equal to 10, while the G and D columns exhibit numbers respectively between 10 and 20 and between 9 and 21. Following the notations in [deH1], the groups of numbers are labelled with an upper case letter for the column and a small case letter for the group. We add the number of notches between parentheses, and get Ma(4) to Mh(7), Da(11) to Dd(9), and Ga(11) to Gd(19).

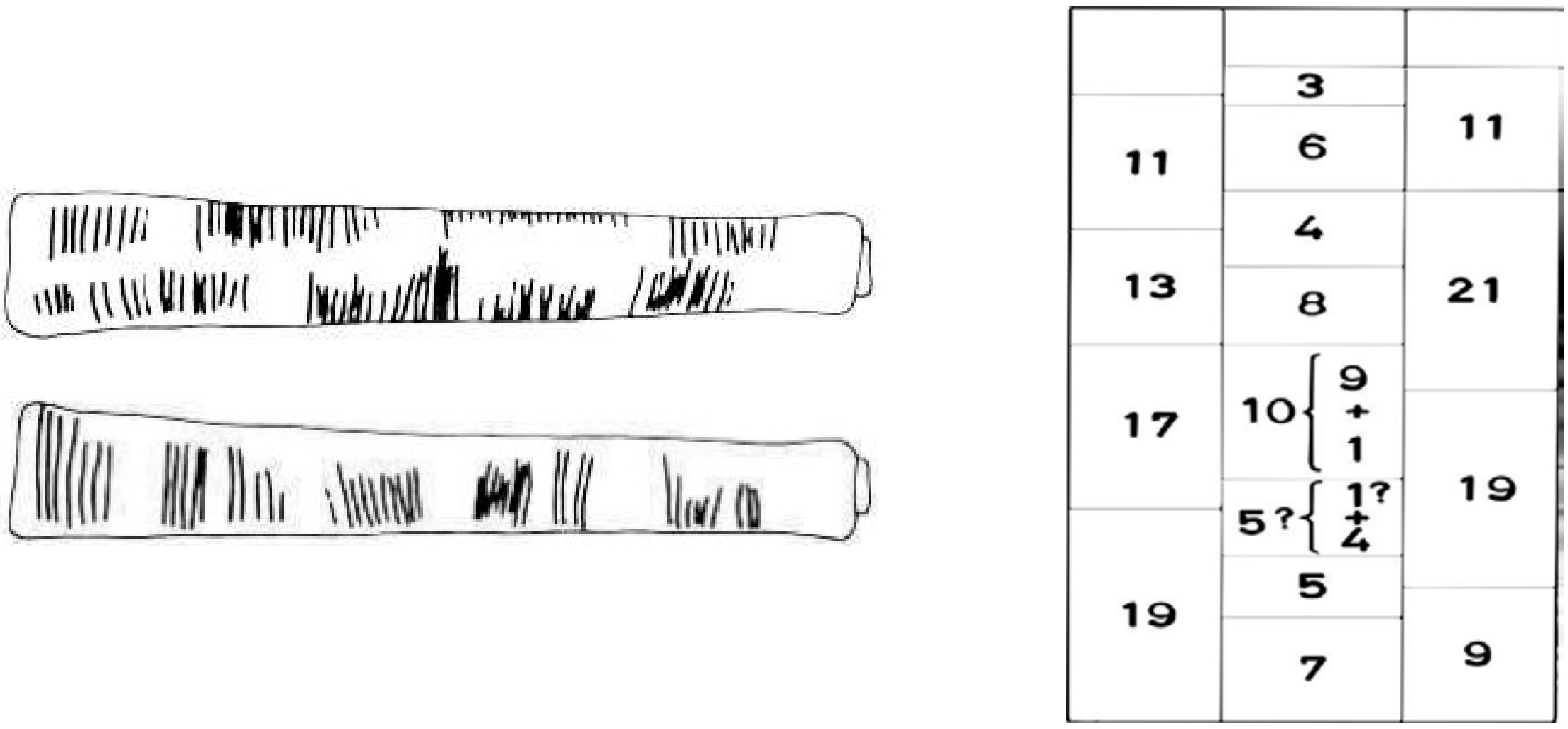

Figure 2: This schematised representation of the notches on the Ishango bone is used throughout the paper.

Early interpretations [deH1, deH2] suggest an ancient knowledge of simple arithmetic. Indeed, the first four groups of the M column suggest duplication by 2while the G column shows prime numbers between 10 and 20. The third column, D, indicates numbers 10 and 20 plus and minus one. In addition, simple operations (addition, subtraction, and multiplication) indicate a wealth of other relations (see e.g. Table V in [deH1]). The most striking is that all numbers in columns G and D add up to 60, while the sum of the numbers in the M column is 48. Since this is respectively 5 x 12 and 4 x 12, the bone discoverer suspected here some simple arithmetic, using the base 2 (for duplication) and the base 10. His major justifications were the central position in the M column, and the display of 10±1 and 20±1 in the D column, along with an early knowledge of prime numbers (G column).
A. Marshack [Mar] suggested a different usage for the bone as a lunar calendar. This surprising alternate interpretation should not be ruled out, especially because of the very convincing circumstantial that can be provided. Indeed, present day African civilisations do use bones, strings and other devices as calendars (see [Lag]). The lack of such evidence was one of the shortcomings in de Heinzelin's arithmetic game interpretation. For example, no awareness of the notion of prime numbers has been discovered before the classical Greek period. The new interpretation that is proposed here combines broadly accepted circumstantial evidence with a straightforward mathematical insight.

## 2. Description of the middle M column

As the middle M column is central to the understanding of the numbering system and the arithmetic displayed on the Ishango bone, a closer examination is in order to try to decipher the way the ancient Ishango people were counting. Furthermore, the operation of replacing the group of notches by the number of notches, although inferring tantalising interpretation, is quite reducing as information is lost mainly about the notch lengths and orientations. It is proposed to re-examine the way these notches are ordered in the central M column.
Table 1 shows the approximate length in mm of each notch and the vertical distance between the eight groups (i.e. the distance along an average top-bottom direction counted from the right tips of notches). Approximate orientations of the different groups are also indicated (with horizontal referring to a direction perpendicular to the average top-bottom direction). Although all distances are approximate (say within ± 1 mm), interesting deductions can be made, when reading from top to bottom (see Figure 2).

TABLE 1
APPROXIMATE ORIENTATION AND DIMENSIONS OF M COLUMN NOTCHES

| Group | Orientation | V | L | Remark on notch length and orientation |
|---|---|---|---|---|
| Ma(3) | horizontal | 20* | 8 | |
| | | | 8 | |
| | | | 8 | |
| Mb(6) | horizontal | 5 | 7 | |
| | | | 4 | |
| | | | 4 | |
| | | | 4 | |
| | | | 7 | |
| | | | 13 | left tip slightly curved downward |
| MC(4) | upward (≈10°) | 18 | 16 | |
| | | | 15 | (or 4+12=16) ; slightly curved upward |
| | | | 14 | |
| | | | 13 | (or 5+8=13) ; slightly curved upward |
| Md(8) | upward (≈10°) | 6 | 13 | (or 7+6=13) ; slightly curved horizontally at middle |
| | | | 12 | |
| | | | 11 | |
| | | | 7 | |
| | | | 7 | |
| | | | 9 | |
| | | | 9 | |
| | | | 9 | |
| Me(9) [Me(10)] | downward (≈10°) | 11 | 11 | |
| | | | 10 | |
| | | | 10 | (or 8+4=12) ; left tip slightly curved horizontally |
| | | | 9 | |
| | | | 9 | |
| | | | 11 | |
| | | | 10 | (or 3+7=10) ; left tip slightly curved horizontally |
| | | | 10 | |
| | | | 13 | |
| | | | [7] | [or 3+2] ; interrupted notch |
| Mf(5) | horizontal | 16 [13] | 3 | notch slightly upward |
| | | | 7 | |
| | | | 9 | |
| | | | 16 | |
| | | | 15 | |
| Mg(5) | horizontal | 5 | 14 | |
| | | | 18 | |
| | | | 16 | |
| | | | 15 | |
| | | | 15 | |
| Mh(7) | horizontal | 10 | 17 | |
| | | | 18 | |
| | | | 18 | |
| | | | 23 | |
| | | | 23 | |
| | | | 23 | |
| | | $5^{+}$ | 22 | (or 17+5=22) ; left tip curved upward |

V: vertical distance (mm) between groups; $^{*}$ from bone top edge ; $^{+}$ to bone bottom edge.
L : Approximate notch length (mm) from tip to tip (or long notch).

The first group Ma(3) is at 20 mm from the upper bone edge; the three notches have the same length (8 mm), are equally spaced and quite parallel in an approximate horizontal direction. The second group Mb(6) is separated vertically from Ma(3) by 5 mm; the six notches have different lengths and are arranged such as the 2nd, 3rd and 4th have a same length (4 mm), with longer notches on each side (1st and 5th, length 7 mm), while the last notch is longer (13 mm); the vertical distances between the 2nd, 3rd and 4th notches (counted from the right tips) are smaller than between the 1st and 2nd, the 4th and 5th, and the 5th and 6th; all notches are approximately horizontal, except for the last one having the left tip slightly curved downward.

The third group MC(4) is separated from Mb(6) by 18 mm; the four notches have approximately the same length (between 13 and 16 mm), are equally spaced (except for the right part of the 2nd notch), and approximately parallel in an upward inclined (right to left) position of approximately 10°, with the 2nd and 4th notches slightly curved upward; the 2nd notch shows a dent at 4 mm from the right tip, most likely due to a workmanship mistake in the notch execution.
The fourth group Md(8) is separated from MC(4) by 6 mm; the eight notches have different lengths and are arranged such as the middle two notches (4th and 5th) are shorter (7 mm), with on either side two subgroups of three longer notches, of approximately equal length in each subgroup, the 1st, 2nd and 3rd notches of lengths between 11 and 13 mm, while the 6th, 7th and 8th notches are 9 mm long; the distances between the three subgroups are visibly larger than the distances between notches within each subgroup; all notches are approximately parallel in an upward inclined direction of approximately 10°, except the 1st slightly curved horizontally at the middle.
The fifth group Me(9/10) is separated from Md(8) by 11 mm; it is not obvious if this fifth group is composed of nine or ten notches; the 1st to 8th notches have lengths between 9 and 11 mm and gives an overall impression of being of approximately equal length; the 9th notch is longer (13 mm); the 10th notch is made of two parts of lengths 3 and 2 mm, separated by a 2 mm gap, yielding a total length of 7 mm between the two part extreme tips; the vertical distances between the 8th and 9th and between the 9th and 10th notches are slightly larger (approximately 3 mm) than among the first 8 notches (approximately 2 mm); the overall orientation is approximately inclined at 10° in a downward (right to left) direction, with the 3rd and 7th notches having their left tips slightly curved horizontally; remark that the 10th notch and the 1st notch of the next group are in a bone area whose surface is damaged.
The sixth group Mf(5) is separated from Me by 16 mm if the Me 10th notch is ignored or by 13 mm is the Me 10th notch is counted; the five notches have different lengths:3 mm for the 1st, a similar length (7 and 9 mm) for the 2nd and 3rd, longer with nearly the same length (16 and 15 mm) for the 4th and 5th notches; the overall orientation is approximately horizontal, with the 1st notch slightly directed upwardly.
The seventh group Mg(5) is separated from Mf(5) by 5 mm; the five notches are arranged in two subgroups, the first composed of the 1st and 2nd notches, the smaller and the longer in this group (14 and 18 mm), the second subgroup is made of the 3rd, 4th and 5th notches of nearly the same length (16 and twice 15 mm); the overall orientation is approximately horizontal, except for the 1st notch, slightly inclined downward (at approximately 3°) and closer to the 2nd notch than the other notches.
The eighth group Mh(7) is separated from Mg(5) by 10 mm; the seven notches are clearly separated in two subgroups of three and four notches, the first three being nearly equal (17 and twice 18 mm), the last four being longer and also nearly equal (three times 23 and 22 mm); the overall orientation is approximately horizontal except for the last notch whose left tip is curved upward; the last notch is at roughly 5 mm from the bone lower edge.

### 3. The M column: some deductions

What can be deduced from this description?
Firstly, that it clearly appears that certain groups have to be considered together, in view of the different vertical separation of the groups (5 or 6 mm and 10 mm or more); the Ma(3) and Mb(6), the MC(4) and Md(8), and the Mf(5) and Mg(5) groups have to be looked at together, while the Me(9/10) and Mh(7) groups could be considered independently (this was already recognised by the bone discover). This is emphasised by the different inclinations of groups that have to be considered together: Ma(3) and Mb(6) are approximately horizontal, MC(4) and Md(8) are inclined upward at approximately 10°, Me(9/10) is inclined downward at approximately 10°, Mf(5), Mg(5) and Mh(7) are approximately horizontal.

Secondly, the first four groups suggest an evident knowledge of duplication, as Mb(6) and Md(8) have a number of notches double of those of Ma(3) and MC(4), which was already noted [deH1]. However, what seems remarkable is that the duplication process is not conducted as just repeating twice the number of notches of an initial group, but more like rearranging a new group in a particular way. Consider Mb(6): start with a set of three notches of the same length, like in Ma(3), although shorter; add two longer notches, one on either side; and add a last even longer notch on only one side to obtain finally six notches, instead of taking two sets of three notches of similar lengths. For the duplication of MC(4), start with a subgroup of two notches in the middle and add on either sides a subgroup of three notches longer than the first two but of similar length within each subgroup; or the other way around, start from two subgroups of three notches of similar lengths and add in the middle two shorter notches. Clearly, the duplication process is not really an operation of making an additional copy added to the original, nor a straightforward multiplication by 2. It is more a reconstruction of an entire group having the required number of notches double than the initial group but arranged differently, each time involving the number 3 and paired subgroups on either side of a central subgroup, that can be schematised by

| | |
|---|---|
| Ma(3): | 3s |
| Mb(6): | 1m + 3s + 1m + 1L |
| MC(4): | 4L |
| Md(8): | 3L + 2s + 3m |

where s, m and L stand for small, medium and long, that can have different length values in different groups.

Thirdly, two interpretations can be proposed for the Me group, depending if one includes the 10th interrupted notch or not, seen either as an initially complete 10th notch damaged with time or as an artefact.
Considering the latter case, the ninth notch is longer than the first eight notches. Based on how their lengths are interpreted, these eight notches can be seen as four notches repeating twice MC(4) (although inclined differently). Or else, they can be explained like Md(8), because of the two slightly shorter notches (9 mm) in the middle and two sets of three slightly longer notches (10 to 11 mm) on either side. Finally, Me(10) could simply be a single block of eight notches since their lengths are all almost equal and most likely within the execution precision allowed by ancient workmanship. In any case, these three interpretations yield a total number of 8, to which are added the longer (13 mm) and slightly separated 9th notch. It suggests how the number 9 is formed by addition of 1 to either 2 x 4, or to 3+2+3, or to 8. In any case this number 9 is not obtained by the process of triplication, i.e. multiplying by 3, but more by duplicating MC(4) or by extending the processes of MC(4) to Md(8), and adding the required notch to make 9.
We now assume the former case of an initially complete notch. Its vertical separation from the 9th notch leads to the interpretation that it forms a part of an initial subgroup of two notches with the 9th notch or that it is an additional single notch, forming the number 10 by either adding 2 to 8 or adding 1 to 9. This latter hypothesis is the one retained by the discoverer of the bone.
Both hypotheses can be schematised as

Me(9): 4m+4m+1L or 3m+2s+3m+1L or 8m+1L
Me(10): 4m+4m+2L(?) or 3m+2s+3m+2L(?) or 8m+2L(?)
or 4m+4m+1L+1L(?) or 3m+2s+3m+1L+1L(?) or 8m+1L+1L(?)

In any case, considering this group Me of 10 notches (if it ever includes ten notches) as central and as representative of the counting base seems not adequate. Another appropriate approach is indicated further.

Fourthly, the two groups Mf(5) and Mg(5) being close together also indicate that they have to be considered together, showing two different ways of obtaining the number 5. The first notch of Mf(5) is much shorter than any of the other in this M column, and most likely was initially longer and was damaged with time. In fact, all the area covering the Me(9/10) last notch and the Mf(5) first three notches show signs of damage. It is hypothesised that initially these three first notches were of approximately the same length, forming a first subgroup of three notches, to which is added a second subgroup of two longer notches to form the number 5.
The Mg(5) group exhibits a separation as well, with the last three notches forming clearly a subgroup with similar lengths (15 to 16 mm) and all approximately horizontal; but again three interpretations are possible for the other two notches. Reasoning on their lengths, one can see them as being two single separate notches to be added to the subgroup of three notches; however, this interpretation is flawed by the closeness of the 1st notch to the 2nd. Their vertical separation of 2 mm between the right tips of the first two notches and between the last three notches, and of 3 mm between the 2nd and 3rd, suggests they are a subgroup of two notches with the first one may be inadequately executed. Reasoning on their orientations and disregarding their lengths (being all close between 14 and 18 mm) and their vertical separation, the 2nd notch seems to belong to the second subgroup to form a subgroup of four notches to which is added the 1st one to form the number 5.
These addition processes in Mf(5) and Mg(5) can be schematised as

Mf(5): 3m (?) + 2L
Mg(5): 1m + 1L + 3m or 2L (?) + 3m or 1m + 4L (?)

Recall that the indication s, m, L are somehow subjective and are different for each group.

Fifthly, the last group Mh(7) shows clearly a separation in two subgroups of three notches and four longer notches, indicating the number 7 as the sum of 3 and 4:

Mh(7): 3m + 4L

Sixthly, a representation of the first two natural numbers 1 and 2 is nowhere to be seen on the bone markings. All other numbers up to 9 (or 10) are represented but surprisingly the first two are missing, although used in the basic calculation operations as addition and duplication.

## 4. The bone of Ishango, another interpretation

All the above considerations lead to the following new interpretation of this M column reading.
The group Me(9/10) occupies a central position in the M column, but it is not the central base on which the arithmetic shown in this column is built. If 10 would have been a natural base on which the arithmetic displayed on the bone would have been built, why is it presented as 9 + 1? This has no logical correspondence to any anatomical feature. One can understand 10 as representing the ten fingers of the two hands, or twice the 4 main fingers plus two thumbs, but not 9 fingers plus one separate finger.
Instead, one should consider the central role played by the two numbers 3 and 4 in all the arithmetical processes considered above. In addition, among the eight groups, the groups Ma(3) and MC(4) are displaying the most regular arrangement with three horizontal notches of equal length and four notches of approximately equal length, inclined upward and parallel to each other (except for the dent in the MC(4) 2nd notch).

Based on the vertical separation of the different groups, it is proposed to re-arrange them in four larger groups or families associated with the simple smaller number exhibited in the first group of each family. Reading the M column from top to bottom, it is proposed:
- to consider the smaller numbers in increasing order, that is 3, 4, 5 and 7;
- to associate to the groups of each of these numbers the groups of other numbers obtained by simple arithmetical operations like duplication and addition of the required number of notches, i.e.:
    - the groups Ma(3) and Mb(6) of 3 and 6 notches (family 1),
    - the groups MC(4), Md(8) and Me(9/10) of 4, 8 and 9 (or 10) notches (family 2),
    - the groups Mf(5) and Mg(5) of 5 notches each (family 3), and
    - the group Mh(7) of 7 notches (family 4).

In the first two families, the associated numbers are obtained simply by either duplication of the first smaller numbers 3 and 4, which are the bases for these two families, and/or successive addition of 1 or 2. (Duplication here should not be taken in the sense known today but as explained above, as the reconstruction of a new group from an initial group by adding other subgroups on either sides of the central one.)
The third family shows two ways of obtaining the 'composed' number 5, based on addition of 1 or 2 to either of the bases 3 and 4. The fourth family shows how to obtain the 'composed' number 7 by adding the two bases 3 and 4.
This concept of regrouping in families is supported also by the different vertical separations: 18 mm between families 1 and 2, 16 mm (or 13 mm) between families 2 and 3, and 10 mm between families 3 and 4 (this distance is smaller most probably because of the proximity of the bone lower edge). This regrouping concept is not essential (and probably did not exist in the mind of the notch maker), but it puts in evidence the fundamental role of the two bases 3 and 4. These two numbers are known since long to play a central role in ancient and present African tribes (see below). The numbers 3 and 4 could have formed the base of the arithmetic system used by the ancient Ishango people for operations on small numbers and that the derived base 12 was used for larger numbers.

## 5. How to account for the G and D columns?

The G and D columns bear four groups of notches each, clearly separated and their interpretation seems quite straightforward. The G column has respectively from top to bottom 11, 13, 17 and 19 notches, while the D column has 11, 21, 19, 9 notches. The interpretation given by the bone discoverer and repeated later by other authors (see for instance [Jos], [Nel] and [Zas]) is to see in the G column the four prime numbers between 10 and 20, suggesting a knowledge of division and of the particular characteristics of prime numbers, and in the D column the four numbers 10±1 and 20±1.
However, several observations do not fit with this simple scheme.

First, the G column shows the four numbers in an increasing value order from top to bottom: 11, 13, 17, 19. If the above idea of the regrouping in four families in the M column is accepted, they are also arranged in an increasing order: 3 (and 6), 4 (and 8 and 9 or 10), 5 (and 5), and 7. Therefore, how to explain the lack of order in the presentation of the 10±1 and 20±1 numbers? Why not 9, 11, 19, 21 instead of 11, 21, 19, 9? Even reading from bottom to top (9, 19, 21, 11) does not show any apparent rationale. On the other hand, if there is a hidden rationale in the sequence 11, 21, 19, 9 in the D column, why does it not appear also in the G column, like e.g. 13, 19, 17, 11?

Secondly, to see prime numbers between 10 and 20 in the G column suggests the knowledge of particular characteristics of these numbers. This however necessitates arithmetical knowledge more advanced than the one displayed on the two other columns, i.e. basic arithmetic operations of addition and multiplication by 2 (or duplication). There is no evidence of knowledge of other operations, like subtraction or division or even multiplication by factors other than 2 (recall that 9 in Me(9/10) is not formed by triplication but probably by duplication and addition of one unity). Furthermore, 3, 5 and 7 are not singled out as being also prime numbers, not mentioning the single even prime 2. The number 2 is not even appearing in the markings of the bone. On the other hand, other evolved though simple concepts like squaring or powers of 2 are not displayed. This is not a prerequisite to the knowledge of prime numbers but it seems easier to display on a counting tool than immediately jump to display prime numbers and furthermore, only those between 10 and 20. Therefore, it appears questionable that the ancient people of Ishango mastered the knowledge and the particularities of the mathematics of prime numbers.

Thirdly, what is the reason for displaying the numbers 10±1 and 20±1 in this disordered manner in the D column? Why the operations of addition and subtraction of 1 although there is no other obvious example of subtraction to be found in the bone markings? The accepted explanation of linking these numbers to a base 10 stresses the influence of the actual modern use of the base 10 more than intrinsic relations to the simple arithmetic displayed in the M column. However, these seem more significant in view of the use of the bases 3 and 4. If these operations of addition and subtraction of 1 from a base are genuinely correct, it would make more sense to see the ordered series of numbers in the G column as displaying the numbers 12±1 and 18±1 from the derived base 12 and the number 18, one and half times the derived base 12.

Fourthly, it is however undeniable that there is an intention behind these markings. The sums of the G and D column notches are equal to each other and furthermore equal to 60. This reinforces the proposed idea of the usage of the derived base 12 as the counting base for large numbers.

## 6. Operations within columns G and D.

To account for these two columns, it is proposed to see a continuation of the arithmetic lesson started on the M column instead of two series of somehow special numbers. The numbers could display further results of simple operations from the basic numbers of the M column, i.e. the bases 3 and 4. The other 'composed' numbers 5 and 7, and the associated numbers 6, 8 and 9 (or 10) could be involved too.

Yet, which operations? There are many combinations possible from these numbers to obtain those displayed in the G and D columns (see e.g. Table V, page 68 of [deH1]). Again, one should be guided first by simplicity, i.e. by trying to find simple relations between the existing numbers, second, by the geometrical arrangement on the bone, and third, by the different patterns of the series of notches. For example, the number 19 shown in the last Gd(19) group is shown as 14 notches separated from 5 notches by 4 or 5 mm, while the vertical separation between other adjacent notches is approximately 2 to 3 mm.

The question why the D column is not displaying a series of notch numbers in an increasing value order should be addressed first. It is proposed to see in this fact not an intended vertical series of numbers but more the results of simple operations carried out on horizontally adjacent numbers. Like in slide rules (used before the apparition of pocket calculators), the notch numbers in the D and G columns could correspond to simple operations involving adjacent numbers in the M column. They are located in approximately the same vertical range as the results displayed in the G and D columns.

Let us start with the central numbers in the D columns, 21 and 19.
The Db(21) group of 21 short notches covers a vertical range of 40 mm, encompassing the MC(4) and Md(8) groups and is close to the Me group of 9 (or 10) notches. The sum of 4 and 8 gives 12, the derived base, to with is added 9 (or 10), yielding 21 (or 22).
The Dc(19) group covers also a vertical range of 40 mm, encompassing the Me(9/10) and Mf(5) groups and is close to the Mg(5) group. The sum of their notch numbers yields 19 (or 20).
Now we turn to the two central numbers of the G column, 13 and 17.
The Gc(17) group covers a vertical range of 33 mm, encompassing the Me(9/10) group and is just under the last notch of the Md(8) group, yielding a sum of 17 (or 18).
The Gb(13) group covers a vertical range of 24 mm, encompassing the MC(4) group and partially the Md(8) group. Their sum gives 12. Another combination is proposed, with the two groups Ma(3) and Mb(6), adjacent to the vertical range of the Gb(13) group. The sum of the first three groups Ma(3), Mb(6) and MC(4) yields then 13.

One has now the beginning of a pattern. Taking consecutive sums of two or three consecutive numbers in the M column starting from the top, their results are displayed in the same approximate vertical range in the middle of the G and D columns:

| M | | G | D |
|---|---|---|---|
| 3 + 6 + 4 | ==> | 13 | |
| 4 + 8 + 9 | ==> | | 21 |
| 8 + 9 | ==> | 17 | |
| 9 + 5 + 5 | ==> | | 19 |

In this case, we have not considered the tenth notch of the fifth group of the M column, assuming it being an artefact.
Pursuing in the same direction and considering that this tenth notch is genuine, the other possibility would give

| M | | | G | D |
|---|---|---|---|---|
| 4 + 8 | *+1* | ==> | 13 | |
| 4 + 8 + 10 | *- 1* | ==> | | 21 |
| 8 + 10 | *- 1* | ==> | 17 | |
| 10 + 5 + 5 | *- 1* | ==> | | 19 |

This second possibility is less elegant (and less simple) as it introduces additionally one addition and three subtractions of 1. The first possibility abides by the guideline rule of simplicity where no additional operation of ± 1 is introduced, but for which the tenth notch of the Me group is not considered. In that case, the partitioning and horizontal equivalence between groups of numbers in the three columns is not respected in the first case (3+6+4 = 13). On the other hand the second possibility has the disadvantage of introducing the additional operations of ±1, but respects the horizontal equivalence between the different columns and considers the tenth notch of the Me group.

Let us see now the four other numbers in the G and D columns, 11 (twice), 9 and 19.
The two first groups Ga(11) and Da(11) are rather different. Although their vertical ranges are about 22 mm, the 11 notches of Da(11) are all parallel and horizontal, while those in Ga(11) have upward inclinations (right to left) between approximately 13 and 24°. The first and last notches of Ga(11) are about 8 to 9 mm below respectively the first and last notches of Da(11). The first notch of Da(11) is approximately 4 mm above the first notch of Ma(3), while the last notch of Da(11) is about at the same level of the last notch of Mb(6). The first notch of Ga(11) ends (or starts) just under the level of the last notch of Ma(3), while the last notch of Ga(11) starts (or ends) approximately 7 mm under the level of the last notch of Mb(6). The notch lengths of Da(11) vary between 4 and 7 mm. The vertical position of Da(11) in correspondence with the first two groups of the M column suggests clearly that 11 should be the result of some operation conducted on the values 3 and 6 of Ma(3) and Mb(6). However, their addition gives 9 and not 11. A closer look at the vertical position of Da(11) reveals that the last 9 notches (3rd to 11th) of this group fall exactly between horizontal lines passing through the 1st notch of Ma(3) and the last notch of Mb(6). This suggests that the result of the operation of adding 3 to 6 in the M column is exactly 9 in the D column first group, to which is added two notches (1st and 2nd) on top of the D column. Therefore, this group of 11 notches is displaying the result (9) of the addition of the first two numbers (3 and 6) in the M column, plus 2. Why the additional adding of 2? No reason can be proposed but it seems that the relative positions of the notches of these three groups are not a coincidence and reflect an unknown intention.
Looking closely at the other group of 11 notches Ga(11), the 1st notch is below the level of the last notch of Ma(3). Furthermore, these 11 notches are formed by three subgroups: the first six notches are approximately parallel in an upward inclined (right to left) position at approximately 13°; the next four notches (7th to 10th) are also approximately parallel with a higher upward inclination of approximately 23° (the difference in inclinations of approximately 10° of these two subgroups recalls the inclination of about 10° of MC(4) with respect to the approximate horizontal direction of Ma(3) and Mb(4)); the last 11th notch of Ga(11) is vertically separated by approximately 3 mm while the distance among the first ten notches are closer to 2 mm. The lengths of the different notches tell us also that these subgroups are somehow differently composed. Schematically, one has

Ga(11): (2m+1L+2m+1L) + 4L + 1L

(recall that m and L can have different values in different groups). So again we have an indication that the values of Mb(6) and MC(4) should be added (6+4), to which is further added 1 notch.

Moving on to the last two groups of the D and G columns, having respectively 9 and 19 notches, Dd(9) has a vertical range of 20 mm, encompassing Mh(7) and under the level of Mg(5). Dd(9) has its first two notches slightly shorter (6 and 7 mm) and inclined upward (right to left) while the other 7 notches are longer (8 to 10 mm) and more or less parallel and horizontal. This subgroup of 7 notches is made of three central (5th to 7th) long notches, with slightly shorter (4th and 8th) notches on both sides and two longer (2nd and 9th) notches at the outside on both sides. This can be represented schematically by

Dd(9): 2s + (1L+1m+3L+1m+1L)

This overall representation suggest that the number 2 is added to the number 7 of Mh(7), like in Da(11) of 11 (2+9) notches.
The G column last group Gd(19) has the longest vertical range of 44 mm and encompasses the last three M column groups Mf(5), Mg(5) and Mh(7). The first 9 notches of Gd(19) are approximately parallel, the next 5 notches slightly curved upward and the last 5 notches slightly shorter and approximately horizontal and parallel again. This Gd(19) group is clearly formed of several subgroups: the first includes the first 2 notches, the second subgroup is made of the next 7 (3rd to 9th) notches; these two subgroups are separated vertically by 3 mm, while the average vertical separation within the two subgroups is approximately 2 mm; the third subgroup, separated by approximately 3 mm from the second one, has its 5 notches curved upward and is made of 3 notches (10th to 12th) slightly separated from the next 2 notches (13th and 14th); the last subgroup, separated by 4 mm from the third subgroup, is made of 5 notches with two shorter (15th and 19th) notches on either sides of three central notches (16th to 18th). Schematically, it gives

Gd(19): 2m + ((7L)+(3L+2m)+(1s+3m+1s))

showing that the number 2 is added to the numbers 7, 5 and 5 in the reversed order of the M column last three groups. The reason for this reversed presentation can not be understood at the first glance but an explanation is proposed further. Interestingly, the two numbers 5 are represented here differently (3+2 and 1+3+1) than in Mf(5) (1+2+2 or (3(?)+2)) and Mg(5) (1+4 (or 2+3)).

## 7. The G and D columns, key to the base 12

If the 10th notch of Me is genuine, one has the following eight addition relations involving three extra additions of 2, two extra additions of 1 and three extra subtractions of 1.

| | M | | | G | D |
|---|---|---|---|---|---|
| *2+* | 3 + 6 | | ==> | | 11 |
| | 6 + 4 | *+1* | ==> | 11 | |
| | 4 + 8 | *+1* | ==> | 13 | |
| | 4 + 8 + 10 | *- 1* | ==> | | 21 |
| | 8 + 10 | *- 1* | ==> | 17 | |
| | 10 + 5 + 5 | *- 1* | ==> | | 19 |
| *2+* | 5 + 5 + 7 | | ==> | 19 | |
| *2+* | 7 | | ==> | | 9 |

For the three operations where the adding of 2 is involved, the two additional notches are above the one, two and three groups of notches of the M column groups. For the two operations where the adding of 1 is involved, the additional notches are below the two groups of the M column groups. Furthermore, these extra additions of 2 and 1 are placed in a symmetric manner at the top and bottom of both G and D columns, while the three subtractions of 1 are involving the middle parts of both columns.

It is interesting to present the M column differently with the different notch groups corresponding to those of the M column under each other. In a descending order (represented here from left to right) and adding (∑) the numbers of the same group appearing several times, this yields

| | M | | | | | | | | | | | G | D |
|---|---|---|---|---|---|---|---|---|---|---|---|---|---|
| *2+* | 3 + | 6 | | | | | | | | | ==> | | 11 |
| | | 6 + | 4 | | | | | | *+1* | | ==> | 11 | |
| | | | 4 + | 8 | | | | | *+1* | | ==> | 13 | |
| | | | 4 + | 8 + | 10 | | | | | *- 1* | ==> | | 21 |
| | | | | 8 + | 10 | | | | | *- 1* | ==> | 17 | |
| | | | | | 10 + | 5 + | 5 | | | *- 1* | ==> | | 19 |
| *2+* | | | | | | 5 + | 5 + | 7 | | | ==> | 19 | |
| *2+* | | | | | | | | 7 | | | ==> | | 9 |
| ∑ = *6* | 3 | 12 | 12 | 24 | 30 | 10 | 10 | 14 | *+2* | *-3* | | 60 | 60 |

It shows that the number 12 appears twice (as twice the group Mb(6) and three times MC(4)), the number 24 appears once (as three times GMC(8)), followed by 30, 10 (twice) and 14. The appearance of twice the derived base 12 and of its multiple 24 is quite significant. The following sum 30 can be seen as a multiple of 6 or as derived from 12 (12 x 2.5) and as half the total sum of each G and D columns. However, the following three sums 10 and 14 are not in line with this approach.
This is possibly why in the last group of the G column Gd(19), the notches are not presented as 2+5+5+7 but as 2+7+5+5. Indeed, putting them in the descending order as presented in the G column (here from left to right), one obtains

| | M | | | | | | | | | | | G | D |
|---|---|---|---|---|---|---|---|---|---|---|---|---|---|
| *2+* | 3 + | 6 | | | | | | | | | ==> | | 11 |
| | | 6 + | 4 | | | | | | *+1* | | ==> | 11 | |
| | | | 4 + | 8 | | | | | *+1* | | ==> | 13 | |
| | | | 4 + | 8 + | 10 | | | | | *- 1* | ==> | | 21 |
| | | | | 8 + | 10 | | | | | *- 1* | ==> | 17 | |
| | | | | | 10 + | 5 + | 5 | | | *- 1* | ==> | | 19 |
| *2+* | | | | | | 7 + | 5 + | 5 | | | ==> | 19 | |
| *2+* | | | | | | | | 7 | | | ==> | | 9 |
| ∑ = *+6* | 3 | 12 | 12 | 24 | 30 | 12 | 10 | 12 | *+2* | *-3* | | 60 | 60 |

This ordering as presented by the marking on the bone in the G column last group yields again the occurrence of twice the derived base 12 (as 5+7). The combination of the sum before last (10) and the two extra added units (*+2*) would yield again the derived base 12, while the first group Ma(3) is nullified by the three extra subtracted units.

## 8. The slide-rule hypothesis (summary).

It is very likely that the ancient people of Ishango made use of the bases 3 and 4 for counting and building further small numbers up to 10. The derived base 12 could have been used for counting larger numbers. It is not impossible that these bases coexisted with other natural counting bases like 10 or 20, although the considerations developed in this paper point toward the bases 3, 4 and 12.
It should be stressed that no evidence of the other arithmetical operations is found in the markings of the bone although addition of simple small numbers was well mastered by the ancient Ishango people. In particular, the multiplication process appears rather primitive as only duplication (or multiplication by 2) could be deduced from the bone markings; knowledge of multiplication by other factors than 2 is not apparent. Furthermore, the hypothesis of mastering the knowledge of prime numbers appears too audacious in view of the basic arithmetic displayed by the bone markings. On the other hand, the proposed hypothesis of considering the bone as an ancient 'slide rule' to display simple addition arithmetic is elegant and fits well with the various notch geometrical patterns. In the additional results displayed on the G and D columns, the adding of extra 2 and 1 and the subtraction of 1 cannot be explained but they show in any case a regular and intended pattern. Finally, the regular occurrence of the derived base 12 in the repetition of the different groups of the M column in the addition operations of the G and D columns show the central role that this base played in the proto-mathematics of the ancient Ishango people.
The precision of this prehistoric 'slide rule' and the regular apparition of the derived base 12 move any mathematician, archaeologist or casual spectator.

## 9. Circumstantial evidence 1: counting methods in Africa.

This slide rule interpretation of the Ishango bone supposes a base 12 number system was used, formed out of combinations of smaller numbers like 3 and 4 together with additions of 1 or 2. One can wonder if such a number system has ever existed. A survey of the numerous varied counting methods in Africa turns this question into a rhetorical one. Zaslavsky wrote a complete book about related subjects, with the transparent title *Africa counts* (see [Zas]). After its publication, many other papers on the subject added even more information. All these publications show that the base ten system is certainly not the only one that was ever preferred. For example, the Congolese Yasayama used a system related to the base 5 (see [Mae]). This can be deduced even today from their number words. They count as follows:

1 = *omoko* 2 = *bafe* 3 = *basasu* 4 = *bane* 5= *lioke*
6 = *lioke lomoko* 7 = *lioke lafe* 8 = *lioke lasasu* 9 = *lioke lane* 10 = *bokama*
11 = *bokama lomoko* 12 = *bokama lafe, etc.*

However, it is not a true quincuncial system, since 25 = 5 x 5 does not seem to play a particular role like 100 = 10 x 10 does in the base ten system. Maybe these number words are only the remains of such a system. The adaptation to the base 10 can have occurred in recent times.
The Congolese Baali system is of a similar type, but now 6 is the base number. It is more noteworthy than the previous since 24 seems to play the role of the number 10 in the decimal system. Indeed, when $24^2$ is reached, a new word is invented, and the construction starts all over again.

1 = *imoti* 2 = *ibale* 3 = *isyau* 4 = *zena* 5 = *boko*
6 = *madia* 7 = *madea neka* (6 + 1?) 8 = *bapibale* (6? + 2)
9 = *bapibale nemoti* (8 + nemoti = 6 + 2 + 1)
10 = *bapibale nibale* (8 + nibale = 8 + 2) 11 = *akomoboko na imoti* (10? + 1)
12 = *komba* 13 = *komba nimoti* 14 = *komba nibale* 15 = *komba nisyau …*
24 = *idingo* 25 = *idingo nemoti …*

| | |
|---|---|
| 36 = *idingo na komba* | 37 = *idingo na komba nemoti …* |
| 48 = *modingo mabale* | 49 = *modingo mabale nemoti …* |
| 576 = *modingo idingo* (= $24^2$) | 577 = *modingo idingo nemoti* (= $24^2 + 1$) … |

A question mark indicated some modified forms that are difficult to explain, but other alterations are easier to understand. The relation between *idingo* and *modingo* is simply one of singular and plural form. In many African languages, prefixes are used instead of suffixes to indicate a plural form (or other relations like diminutives).
Combining number bases, like with 3 and 4 as proposed in the previous sections, has been frequently noticed too. The Nyali from Central-Africa used a mixed system forming numbers through combinations of 4, 6 and 24 = 4 x 6:

| | | | |
|---|---|---|---|
| 1 = *ingane* | 2 = *iwili* | 3 = *iletu* | 4 = *gena* |
| 5 = *boko* | 6 = *madea* | 7 = *mayeneka* | 8 = *bagena* (= plural form of four) |
| 24 = *bwa* | 576 = *mabwabwa* (= $24^2$). | | |

Such combinations are also found with even larger numbers. Linguists who studied the language of another tribe in about the same region, the Ndaaka, discovered the use of bases 10 and 32. The Ndaaka express 10 as *bokuboku*, and 12 as *bokuboku no bepi,* but for 32 there is a special word, *edi*. Now 64 becomes *edibepi* (= 32 x 2) while 1024 is *edidi* (or $32^2$). A number like 1025 is expressed as *edidi negana* or $32^2 + 1$.
Number words can vary according to what has to be counted. For example, if in Burundi one has to do with large quantities of cattle, groups of five have to be used. In this case, the usual word *itandatu* or 3+3 changes in *itano n'umwe,* meaning 5+1. Similarly, *indwi* or 7 changes into *itano n' iwiri* or 5+2, etc. This can come as a surprise to the uninitiated reader, but it should not be more surprising than the changes in words in other languages to indicate for example a male or a female item.
One of the reasons for the great creativity found in the vocabulary for number words can be the various counting gestures. These can still be seen on market or other public places people where they are used to facilitate communication. It is even very likely that the number words are simple reflections or 'translations' of their physical equivalent, but as the Shambaa example will show, this is not necessarily the case.
Zaslavsky pointed out that sometimes stretched fingers designate numbers (like in Western Europe), while sometimes it are the folded fingers. The Soga tribe shows 6 by holding the left forefinger to the closed right hand, while 7 corresponds to adding the middle finger to that left forefinger. The Chagga take the fingers of the right hand, beginning with the little finger, with the whole left hand, to indicate numbers from 6 to 9. The Tete cross the fingers they want to show with the left thumb. Others show a closed fist with the thumb in the middle between the middle and the ring finger to indicate 5, and it is then understood as 2+1+2. Some Maasaï push the top of the middle finger on the thumb and on the top of the forefinger to indicate 3, but if the stretched forefinger rests on the stretched middle finger, it means 4. In Rwanda and Western Tanzania, four corresponds a forefinger against the ring finger of one hand and snapping these fingers onto the middle finger.
A remarkable fact in the denominations and gestures for the numbers from six to nine is that these can be formed by different principles. Sometimes two equal expressions are used, and 6 and 8 are then expressed as 3+3 and 4+4. Other descriptions are frequent for the case of 7 and 9, like 10-3 and 10-1. Still, in some regions compositions based on 5 are preferred: 6 = 5+ 1, 7 = 5+ 2, etc.
An example of different words and gestures for counting from 6 to 9 is the following *Shambaa* - language.

| Number | Spoken | Meaning | Gesture |
|---|---|---|---|
| 6 | *mutandatu* = *ntatu na ntatu* | 3 and 3 | Extending the three outside fingers of each hand, or 3+3. |
| 7 | *mufungate* = *funga ntatu* | take off 3 fingers 7 = 10-3 | Four fingers on the right hand and three on the left, or 4+3. |
| 8 | *munane* = *ne na ne* | 4 and 4 | Four at the right hand, and four at the left hand, or 4+4 |
| 9 | *kenda* | take off one 9 = 10-1 | Five at the right hand, and four at the left or 5+4. |

Finally, it may be noted that such number word constructions were not Africa's exclusivity. The Quevedo of South-America count from one to ten as follows: 1, 2, 3, 2+2, 2+3, 2x3, 1+2x3, 2x4, 2x4+1, 2x4+2.
The construction methods of the number words already provide some circumstantial evidence for present 'base 12 slide rule' interpretation of the notches on the bone. The Baali system showed how 7, 8, 9, 10 and 11 are formed through sums of simpler numbers. The Nyali's 4-6 mixed system illustrated the possibility of a similar 3-4 number construction. Zaslavsky's account on the different gestures showed how number words are expressed physically, and maybe this was precisely what the early Ishango man was doing, consciously or not, when he was carving notches on the bone.

## 10. Circumstantial evidence 2: base 12 counting methods

The slide rule hypothesis affirms more however: the number system constructed through additions and combinations would have involved the number 12 in particular. It happens that one of the many African counting methods similar to the ones given above uses the base 12 in an unmistakable way. In addition, it explains why this duodecimal system very often comes together with references to the sexagesimal base 60, like in the number of minutes in an hour and the number of hours in a day.
The duodecimal base seems as old as the base 10 system. Some five thousand years ago, the Babylonian already divided the zodiacal circle into $360^{o}$, or twelve signs of the zodiac of $30^{o}$ each. In the fourth century AD, Theon of Alexandria pointed to the computational conveniences of using the base 60 and in more recent times, King Carl XII of Sweden even wanted to impose the duodecimal system as the only legal system. Today, remains of that non-decimal system still occur in many Western expressions, such as in the words dozen and gross.
Yet, the origin of the base 12 and of the related base 60 is an often-recurring question, even to non-mathematicians. The usual arithmetic (based on the divisors of 12) and or astronomical explanations (based on the number of moon-months) both are a *posterior*. To realise the advantages of the numerical system, one has to develop such a system first.
Still, a counting technique that considers parts of the fingers to represent the numbers from 1 to 12, is still in use in Egypt, Syria, Turkey, Iraq, Iran and Afghanistan, Pakistan, Indochina, India. The thumb of a hand counts the bones in the fingers of the same hand (See Figure 3). Four fingers, with each three little bones, evidently yield 12 as a counting unit. The thumb itself is the counting tool, and its bones are not considered. Also, each dozen is counted by the fingers of the other hand, including the thumb, and the multiple 5 x 12 = 60 provides an additional indication of the often simultaneous occurrence of the duodecimal and sexagesimal base.

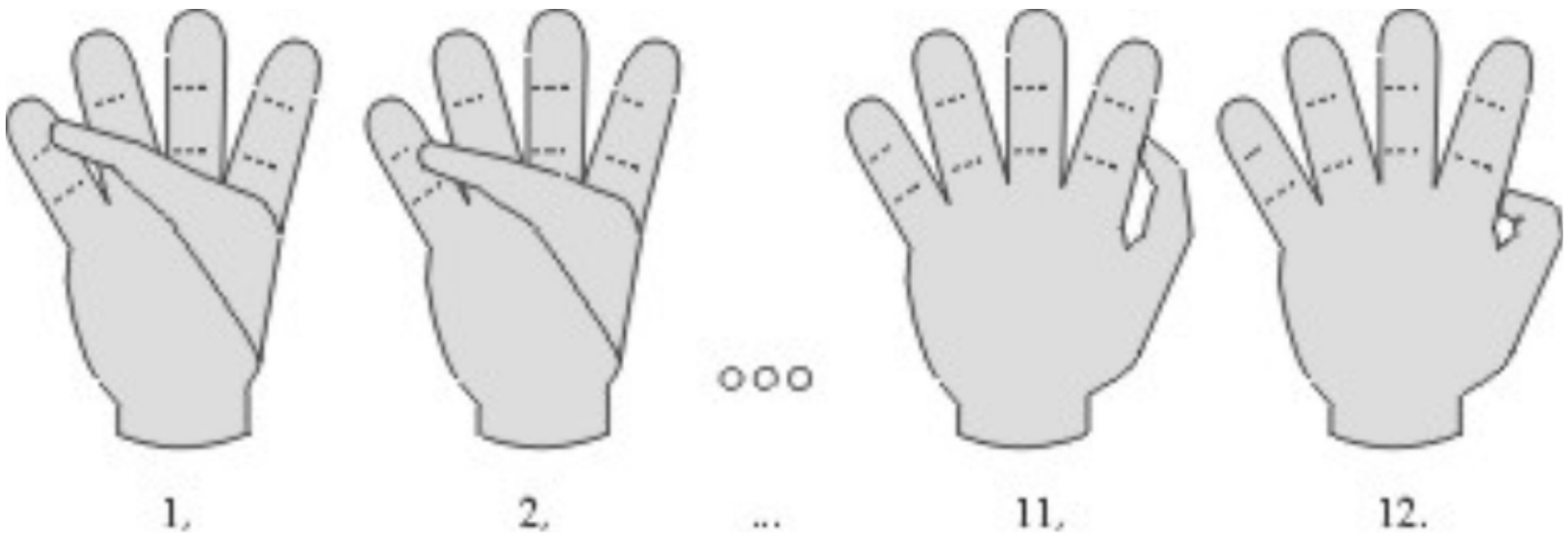

Figure 3: A straightforward explanation for the base 12 uses the thumb as the counting tool, and touches each of the three phalanxes of the four other fingers. his yields a total of 12, and repeating this process for each of the finger of the other hand, the total becomes 48 for the four fingers (fore, middle, ring, and little) or 60 for the complete other hand.

This explanation is not *posterior* like the arithmetical or the astronomical ones. This duodecimal base was indeed a practical one for what these early civilisations wanted to count or to represent. In the matriarchal societies, they could associate the number 1 to the woman, the number 3 to the man, and 4 to the union of woman and man. Or, in after some rather general evolution, they designated the male genitals by the number 3, and the genital symbol of women by 4, making 7 the symbol of their union (see [Nic]). The number 4 seems to have been the most widespread of the mystical numbers. It was established by associations with colours, with social organisation, and with various customs among numerous tribes. The use of six as a mystical or sacred number was less extensively distributed through history and throughout the world than the four-cult, but sometimes a mythology past from quarters cult to a six cult. For example, the four cardinal points (such as North, South, East, West) are simply augmented by the addition of two other points (such as the zenith, above and nadir, below). On the other hand, the counting skills they obtained in this way, allowed them to note that there are 13 (moon)-months in one (solar) year, and not 12 (see [McG]).

This physiological explanation for the duodecimal base is only a hypothesis, but number words as present day tribes in Africa use them, provide further evidence. N. W. Thomas [Tho] reported on such number words in his study of the West-African tribes in the region of the actual Nigeria. Between the rivers Benue and Gurara, which flow into the river Niger more westwards, live the Yasgua, the Koro and the Ham. The Yasgua, Thomas asserted, count as follows:

| | | | |
|---|---|---|---|
| 1 = *unyi* | 2 = *mva* | 3 = *ntad* | 4 = *nna* |
| 5 = *nto* | 6 = *ndshi* | 7 = *tomva* | 8 = *tondad* |
| 9 = *tola* | 10 = *nko* | 11 = *umvi* | 12 = *nsog* |
| 13 = *nsoi (=12+1)* | 14 = *nsoava (=12+2)* | 15 = *nsoatad* | 16 = *nsoana* |
| 17 = *nsoata* | 18 = *nsodso* | 19 = *nsotomva* | 20 = *nsotondad …* |

The Koro do not repeat their word for 12, but another prefix "pl":

| | | | |
|---|---|---|---|
| 1 = *alo* | 2 = *abe* | 3 = *adse* | 4 = *anar* |
| 5 = *azu* | 6 = *avizi* | 7 = *avitar* | 8 = *anu* |
| 9 = *ozakie* | 10 = *ozabe* | 11 = *zoelo* | 12 = *agowizoe* |
| 13 = *plalo (=12+1)* | 14 = *plabe (=12+2)* | 15 = *pladsie* | 16 = *planar* |
| 17 = *planu* | 18 = *plavizi* | 19 = *plavita* | 20 = *plarnu …* |

The Ham tribe proceeds similarly, but again with a different vocabulary. Also, Thomas testified that the base six for numbers from seven to nine was used in Guinea. He reported that 4-6 combinations like 10=6+4 was in use with the Bulanda tribe, while the Bola would have expressed 12 as 6x2 and 24 as 6x4.
This apparently well-known use of base 12 (or base 6?) systems in some areas, surprising to many Western readers, makes the base 12 hypothesis more plausible.

### 11. The Ishango bone as missing link in the base 12 dissemination.

The discoverer of the bone, J. de Heinzelin, also formulated a hypothesis about the impact of the bone. He supposed there could be a relation between the arithmetic operations apparent from the Ishango bone and the dawn of mathematics in ancient Egypt. He documented his assertion by a comparison of harpoon heads found in Ishango and similar objects found in Northern Sudan and Egypt:

> '*The first example of a well-worked-out mathematical table*', says de Heinzelin, '*dates from the dynastic period in Egypt. There are some clues, however, that suggest the existence of cruder systems in predynastic times. Because the Egyptian number system was a basis and a prerequisite for the scientific achievements of classical Greece, and thus for many of the developments in science that followed, it is even possible that the modern world owes a great one of its greatest debts to the people who lived at Ishango. Whether or not this is the case, it is remarkable that the oldest clue to the use of a number system by man dates back to the central Africa of the Mesolithic period. No excavations in Europe have turned up such a hint*'.

The Belgian archaeologist had many followers who amplified his simple hypothesis to incomparable proportions, and even to the hypothesis that there were more than some cultural influences between black Africa and the ancient Egyptians. In 1976 A. Noguera gave his book "*How African Was Egypt?*", the subtitle "*A Comparative Study of Ancient Egyptian and Black African Cultures*", while in 1990 Bernal defended in his '*Black Athena*' the hypothesis that Africa's influence even reached to Greece [Ber1], [Ber2]. There is no scientific consensus about the factual verification of the African-Egypt supposition, which is difficult to substantiate because of the lack if written sources. Up to this present day, the discussion continues and divides the community of archaeologists, Egyptologists and African-Americans in the supporters and opponents. Remarkably enough, this debate is hardly ever heard in Europe, and the subject is not even known to the many Africans in the old continent. In 1995 I. Van Sertima [Ser] published a further defence of Bernal's theorems in his '*Black Athena Revisited*', but barely one year later, the counteroffensive came from M. Lefkowitz in '*Not Out of Africa, How Afrocentrism Became an Excuse to Teach Myth as History*' (see [Lef1] and [Lef2]). The '*Black Athena Debate*' is sometimes very harsh and confused by political statements (see [Row]). Nevertheless, what remains undisputed, is the Ishango bone's title of oldest mathematical artefact.
Yet, mathematics did not necessarily develop in a linear way from Ishango through Egypt and the Middle East to Europe. Maybe the Egyptians had to reinvent the knowledge they had heard about from more southern regions. More generally, sometimes useful skills have been simply forgotten and had to be rediscovered, as Renaissance showed but too well to the Western civilisation. There are sufficient proofs of other arithmetical operations that did develop in one region at a certain time, and not at another region. The discoverer of the bone, who cannot be suspected of not appreciating the marvellous scientific artefact, stated it as follows (see [deH2]):

> *One should not even think about claiming the honour of some universal posthumous certificate in the name of the Ishango man. On the contrary, I personally believe that*

> *most of the civilisations that came afterwards and that have modern affinities, beginning with the final Palaeolithic, the Mesolithic and finally those of the Neolithic stages, those of the agricultural villages and the cities, had to invent several times the same thing, at about the same time, in distinct locations, yet at different periods of time. It is the main idea that imposes itself, for the invention of ceramics, writing, numeration and domestication of animals.*

The Ishango itinerary provides examples of both situations. From Central-Africa to Egypt the message apparently did pass, but there is a part of itinerary proposed by de Heinzelin that seemed to have sunk into oblivion. Indeed, on his map there is a branch from Ishango to West Africa, again obtained through findings of harpoon heads found in Lake Chad and farther Westwards to the region of the upper Niger and even as far as Bamako (see Figure 4). This time, the influence of the mathematical knowledge of the Ishango bone did not seem to have travelled with the harpoon heads, like it was supposed to be the case for its odyssey up north. Actually, this is what one could conclude from the existing literature, but the present 'base 12 slide rule' interpretation of the Ishango bone that puts Thomas' statement in a different light.

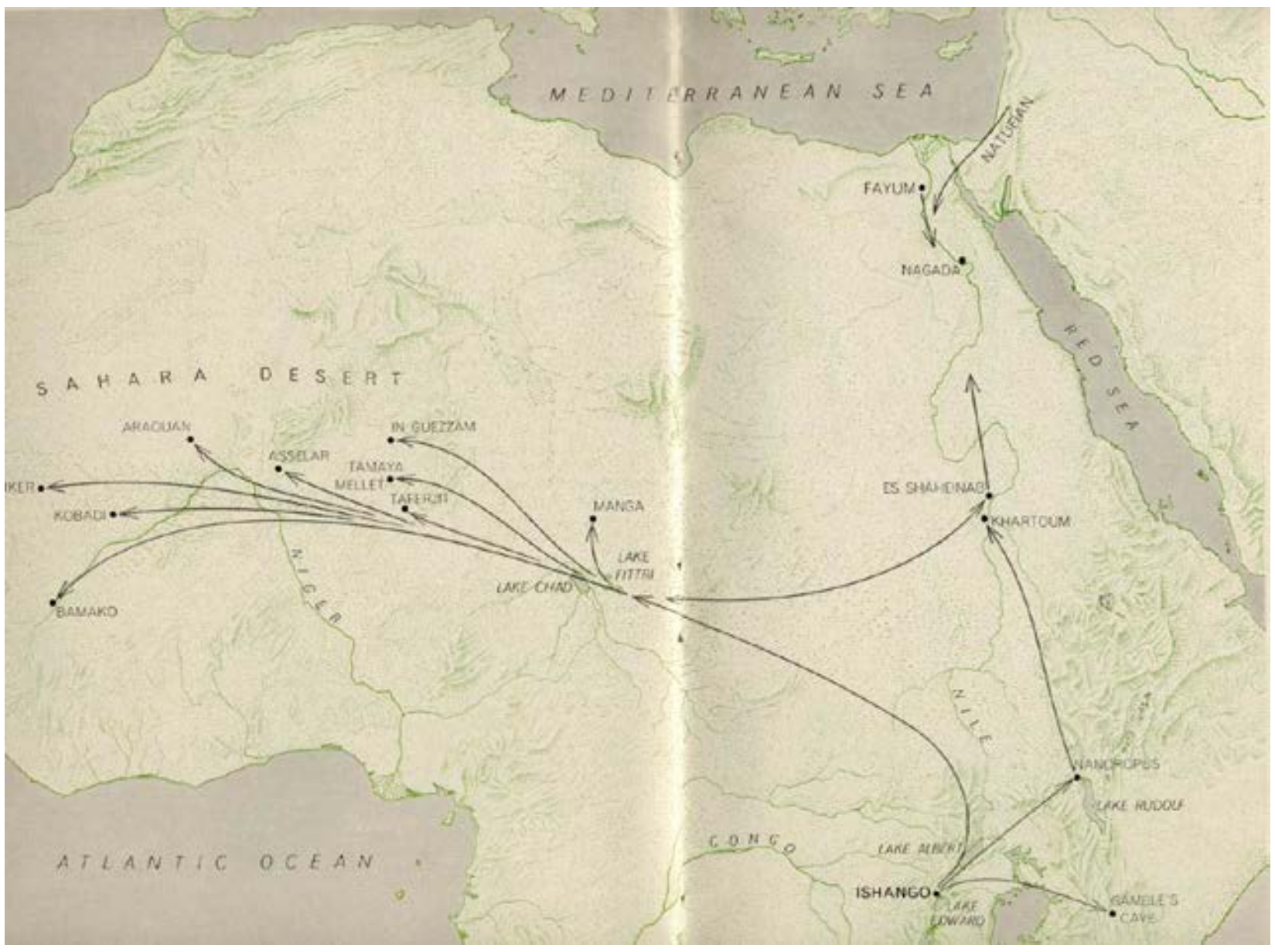

Figure 4: The idea of the influence of non-European cultures on the dawn of science is approved in broad circles of historians of science. Joseph states that the dawn of mathematics is situated in Greece, and that Euclides' creation reached Europe through and thanks to the Arab world [Jos]. Yet, he asserts, Egypt, Mesopotamia, China and India (and the American continent) did obtain some remarkable mathematical results. On this drawing arrows indicate the alternative routes for the development of mathematics. A more far-reaching theory about the influence of Africa on ancient Egypt is not without criticism, but some comparisons, like between the represented Cretan and Nigerian statue, or between the Egyptian and Ugandan traditional dresses, talk for themselves (see [Nog]).

In his paper about West-African tribes, written in 1920, Thomas asserted that the base 12 system had to be rather old, from before the fourteenth century as appeared from the written source "*Travels of Ibu Batuta*". Furthermore, he supposed that in view of European prehistoric discoveries, it should go back two thousand years at least. He knew about Babylon as "*the best known reference to such a base*", but wondered if there was a link between both. He suspected this, since "*as regards burial customs the foreign element is conspicuous*." Although Thomas made rather critical comments about his own findings, he concluded:

*It has frequently been assumed that the duodecimal system, which is in Europe crossed with the decimal system, is a product of Babylonia; how far this view is still accepted I do not know. But it is clear that, even though Egyptian influence in West Africa may be well established, we can hardly accept such a far-reaching theory as Babylonian influence on numbers below 20, which would surely imply both early and close contact, in the absence of other evidence of Asiatic influence in this area.*

*It remains to add that if we find no duodecimal system among any people likely to have been in contact with Nigerian tribes, we must assume an independent origin for the system. If it had been transmitted from Babylonia via Egypt, it must surely have left some traces on its road. For those who believe the duodecimal notation can have been invented once only, it is an interesting problem to bring the Nigerian duodecimal area into relation with Babylonia.*

The Ishango bone could be the missing link Thomas was looking for in 1920. It was found 30 years later, and de Heinzelin made his map about influence from Ishango to Egypt or to West-Africa only in 1962. The archaeologist did this most probably unaware of Thomas' findings, which could be catalogued in more ethnological fields. In 1999, in this actual paper, Pletser interpreted the bone as a base 12 pattern and this again without having read Thomas' work. He actually is a microgravity space science researcher, who was more involved in trying to get the Ishango bone into space (see note 2). He stumbled over his slide-rule interpretation when he was informed about the Ishango bone space project. There seems to be too much coincidence involved in these three research results for them to be sheer accident.

Finally, 80 years ago, Thomas already guessed where he had to look for that missing link. He had pointed out that there is one tribe in Africa with computational practices similar to those he studied in West-Africa. The features of doublings and additions Pletser discovered on the bone are similar to their system of expressing numbers. For instance, 7 is 6+1, while 8 corresponds to 2x4, and 16 to (2x4)x2. The following three numbers are thought of as additives, but 20 is again 10x2. The inventors of this system are the Huku-Walegga, who live in an area Northwest of the Lower Semliki, the same river on which shores the Ishango bone was found.

## Acknowledgements

The author acknowledges the invaluable help of Prof. D. Huylebrouck of the *Sint Lucas* Institute of Brussels for providing the content of sections 9 to 11. Archaeologist A. Hauzeur and Prehistorian I. Jadin of the Royal Museum of Natural Sciences of Belgium in Brussels helped the author and Prof. Huylebrouck at many occasions. In particular, the museum provided the voluminous original report [deH1] of Prof. de Heinzelin and drew the attention to a paper from a century ago about.

## Appeal to the Reader

Like the mathematical knowledge was maybe transferred from Central Africa to Ancient Egypt and then to Europe, the bone followed a somewhat similar odyssey. It was excavated after thousands of years and brought to Europe by its discoverer. This odyssey is not completed since it may fly one day in space. It would be a gesture and a symbol to link early mathematical knowledge and modern technology, mimicking the ellipse in Stanley Kubrick's acclaimed movie "2001: A Space Odyssey" (see [Huy2]).

## References

[Ber1] Bernal M., *Black Athena*, Volume 1: *The Afroasiatic Roots of Classical Civilisation, The Fabrication of Ancient Greece*, Rutgers University 1987.

[Ber2] Bernal M., *Black Athena*, Volume II: *The Archaeological and Documentary Evidence*, Rutgers University 1991.

[deH1] De Heinzelin de Braucourt J., *Exploration du Parc National Albert : Les fouilles d'Ishango*, Fascicule 2, Institut des Parcs Nationaux du Congo Belge, ("Exploration of the National Park Albert : The Ishango excavation", Volume 2, Institute of National Parks of the Belgian Congo), Brussels, 1957.

[deH2] De Heinzelin de Braucourt J., *Ishango*, Scientific American, 206-6, June, 105-116, 1962.

[Huy1] Huylebrouck D., *The bone that began the Space odyssey*, The Mathematical Intelligencer, 18-4, 56-60, 1996.

[Huy2] Huylebrouck D., *Counting on hands in Africa and the origin of the duodecimal system*, Wiskunde en Onderwijs, no 89 Jan-Feb-Mar, 1997.

[Ifr] Ifrah G., *From One to Zero: A Universal History of Number*, New York: Viking, 1985.

[Jos] Joseph G. G., *The Crest of the Peacock: Non-European Roots of Mathematics*, Penguin Books, London, 1992.

[Lag] Lagercrantz S., *Counting by Means of Tally Sticks or Cuts on the Body in Africa*, Anthropos 68, p. 37, 1973.

[Lef1] Lefkowitz Mary R. and Rogers Guy MacLean (eds), *Black Athena Revisited*, University of North Carolina Press 1996.

[Lef2] Lefkowitz Mary R., *Not Out of Africa, How Afrocentrism Became an Excuse to Teach Myth as History*, BasicBooks, 1996.

[Mae] Maes Dr. J., *Note sur les Populations Lalia et Yasayama du Territoire des Dzalia-Boyela*, Congo vol. I-2, p. 172 – 179, 1934.

[Mar] Marshack A., *Roots of Civilisation, The Cognitive beginnings of Man's First Art, Symbol and Notation*, McGraw-Hill Book Company, New York, 1972.

[McG] McGee W. J., *Primitive Numbers*, Bureau of American Ethnology, 1900.

[Moi1] Moiso B. and Litanga N., *Numération Cardinale dans les Langues Bantu du Haut-Zaïre*, Annales Aequatoria 6, 189 – 196, 1985.

[Moi2] Moiso B., *Etude comparée du système de numérotation de 1 à 10 dans quelques Langues non-Bantu du Haut-Zaïre*, Annales Aequatoria 12, 475 – 479, 1991.

[Nic] Nicolas F.-J., *La mine d'or de Poura et le nombre quatre-vingt en Afrique occidentale*, Notes Africaines, Notes Africaines, n° 161, p. 19-21, 1979.

[Nel] Nelson D., G. G. Joseph en J. Williams, *Multicultural Mathematics, Teaching Mathematics from a Global Perspective*, Oxford University Press, 1993.

[Nog] Noguera A., *How African Was Egypt? A Comparative Study of Ancient Egyptian and Black African Cultures*, Vantage Press, New York, 1976.

[Row] Rowe W. F., *School Daze: A Critical Review of the 'African-American Baseline Essays' for Science and Mathematics*, The Skeptical Inquirer, Sept.-Oct, 1995.

[Tho] Thomas N. W., *Duodecimal Base of Numeration*, Man, Nos. 13-14 February 1920.

[Ser] Van Sertima I. (Editor), *Black Athena Egypt Revisited*, London and New Brunswick, Transactions Publications, 2nd Ed. 1995.

[Zas] Zaslavsky C., *Africa Counts*, Lawrence Hill Books, New York, 1973.